\documentclass[12pt]{article}

\usepackage[cp1251]{inputenc}
\usepackage[russian]{babel}
\usepackage{graphicx,comment}
\usepackage{amsmath,amssymb,amsfonts,amsthm,comment}
\newtheorem{theorem}{{\rm Т е о р е м а}}
\newtheorem{lemma}{{\rm Л е м м а}}

\newcommand{\e}{\mathbf{e}}

\begin{document}

\begin{center} {\bf   Голубев Г. К., Крымова Е. А.} \end{center}

\vskip 0.4cm

\begin{center}
{\large\bf ОЦЕНИВАНИЕ  УРОВНЯ ШУМА В ЛИНЕЙНЫХ МОДЕЛЯХ БОЛЬШОЙ РАЗМЕРНОСТИ}
\footnote{Исследование выполнено при частичной финансовой поддержке РФФИ в рамках проекта 15-07-09121
и Deutsche Forschungsgemeinschaft (проект SFB 823, Statistical modelling of nonlinear dynamic processes).}
\end{center}

\bigskip

{\begin{quotation} \small
Рассматривается задача  оценивания уровня шума $\sigma^2$ в гауссовской линейной модели 
$Y=X\beta+\sigma \xi$, где $\xi\in \mathbb{R}^n$ -- стандартный дискретный белый гауссовский шум, а $\beta\in \mathbb{R}^p$
-- неизвестный  мешающий вектор. 
Предполагается, что  $X$ -- известная, плохо обусловленная $n\times p$\,-\,матрица с   $n\ge p$. При этом
 размерность $p$ является большой.
В этой ситуации  вектор $\beta$  оценивается с помощью  спектральной  регуляризации оценки максимального
правдоподобия и  оценка уровня шума  вычисляется  с помощью адаптивного, т.е. основанного на наблюдаемых данных, нормирования   квадратичной ошибки предсказания. Для этой оценки вычисляется скорость ее концентрации   вблизи  псевдо-оценки  $\|\sigma\xi\|^2/n$.     
\end{quotation}}

\vskip 0.7cm

\section
{ Введение и основной результат} 

В настоящей работе рассматривается задача оценивания уровня шума $\sigma^2\in \mathbb{R}^+$  по наблюдениям 
\begin{equation}\label{eq.1}
Y=X\beta +\sigma\xi,
\end{equation}
где $X$ -- известная $n\times p $\,-\,матрица с $n\ge p$, $\beta\in \mathbb{R}^p$ -- неизвестный мешающий вектор,  а $\xi\in\mathbb{R}^n$ -- стандартный дискретный 
белый гауссовский шум, т.е. вектор, компоненты которого являются независимыми гауссовскими случайными величинами с нулевым средним и единичной дисперсией. 

Стандартный подход к оцениванию $\sigma^2$ основан на методе максимального правдоподобия, 
который дает следующую оценку: 
\begin{equation*}
\begin{split}
\widehat{\sigma}^2(Y) =&\arg\max_{\sigma^2>0}\max_{\beta\in \mathbb{R}^p}
\biggl\{-
\frac{\|Y-X\beta\|^2}{2\sigma^2}-\frac{n}{2}\log(\sigma^2)\biggr\}
= \frac{\|Y-X(X^\top X)^{-1}X^\top  Y\|^2}{n}.
\end{split}
\end{equation*}
Статистические свойства этой оценки хорошо известны и могут быть легко 
установлены с помощью метода главных компонент. Пусть $\e_k\in \mathbb{R}^p,\ k=1,\ldots,p$ и $\lambda_1\ge \lambda_2\ge, \ldots,\ge \lambda_p$ -- собственные векторы и собственные числа матрицы $X^\top X$, т. е.
\[
X^\top X \e_k =\lambda_k \e_k,\quad k=1,\ldots,p.
\] 
Определим следующие векторы в $\mathbb{R}^n$:
\[
\e_k^*=\frac{X\e_k}{\sqrt{\lambda_k}}.
\]
Нетрудно проверить, что $\e_k^*,\,  k=1,\ldots,p$, являются ортонормальными векторами в $\mathbb{R}^n$. Дополним эту систему векторов  какими-нибудь другими ортонормальными векторами $\e_k^*, \, k=p+1, \ldots, n$ до полной ортонормальной системы.  

Тогда для линейных статистик $\bar{Y}_k=\langle Y, \e_k^*\rangle$ получим следующее представление:
\begin{equation}\label{eq.2}  
\begin{split}
\bar{Y}_k=&\sqrt{\lambda_k}\bar{\beta}_k+\sigma\xi_k', \quad k=1,\ldots,p,\\
\bar{Y}_k=&\sigma\xi_k',\qquad\quad\qquad  k=p+1,\ldots,n,
\end{split}
\end{equation}
где $\bar{\beta}_k=\langle \beta, \e_k\rangle $, а $\xi'$ --- стандартный белый гауссовский шум. Поскольку $\{\e_k^*, \, k=1,\ldots,n\}$ -- полная система векторов в $\mathbb{R}^n$, то очевидно, что статистические модели \eqref{eq.1} и \eqref{eq.2} являются математически эквивалентными. 
Заметим однако, что эта эквивалентность  уже не бесспорна  с вычислительной точки зрения. Дело в том, что для того, чтобы формально перейти от   \eqref{eq.1} к \eqref{eq.2} нужно вычислить 
спектральное разложение $X^\top X$. Это может быть в принципе достаточно дорогостоящая операция при 
больших $p$ и плохой обусловленности $X^\top X$. Как мы увидим далее, вычисление спектрального разложения часто не является необходимым. Поэтому в этой статье мы будем использовать представление \eqref{eq.2} в основном  для доказательств математических результатов, а строить оценки, используя наблюдения \eqref{eq.1}.

В силу эквивалентности \eqref{eq.1} и \eqref{eq.2} легко проверить, что 
\begin{equation*}
\widehat{\sigma}^2(Y)=\frac{1}{n}\sum_{k=p+1}^n\bar{Y}_k^2=\frac{\sigma^2}{n}\sum_{k=p+1}^n \xi'^2_k.
\end{equation*}
и мы видим, что, во-первых,  оценка максимального правдоподобия $\widehat{\sigma}^2(Y)$ является смещенной. Поэтому, как правило, в  статистике  используется ее несмещенная модификация
\[
\widetilde{\sigma}^2(Y)=\frac{\|Y-X(X^\top X)^{-1}X^\top Y\|^2}{n-p}=\frac{1}{n-p}\sum_{k=p+1}^n\bar{Y}_k^2.
\]   
Во-вторых,  ни $\widehat{\sigma}^2(Y)$, ни $\widetilde{\sigma}^2(Y)$ не используют данные $\bar{Y}_k,\, k=1,\ldots,p$. 
В задачах классической статистики $p\ll n$ и  это
не приводит к существенным потерям качества оценивания. 
В то же время существует достаточно широкий класс статистических задач, в которых   $p\approx n$ и при этом собственные числа $\lambda_k$ достаточно быстро убывают.  В этих задачах мы, естественно, уже не сможем хорошо оценивать $\sigma^2$ 
если не будем использовать  $\bar{Y}_k,\, k=1,\ldots,p$. 

Классическим примером такой задачи является оценивание
$\sigma^2$  в нелинейной регрессионной модели
\begin{equation}\label{eq.3}
Y_i=f(X_i)+\sigma \xi_i, \quad i=1,\ldots,n,
\end{equation}
где $\xi$ -- стандартный белый шум, $X_i\in [0,1]$ -- известные регрессоры, а $f(x), x\in[0,1]$ -- неизвестная мешающая гладкая функция, принадлежащая, например, классу Соболева 
\begin{equation*}
\mathcal{W}_2^m=\biggl\{f: \int_0^1 [f^{(m)}(x)]^2\le L\biggr\};
\end{equation*} 
здесь $f^{(m)}(x)$ -- производная порядка $m$ функции $f(x)$. Более подробно эта модель будет рассмотрена в параграфе \ref{Re}.

Другим хорошо известным  примером являются  обратные задачи, возникающие при дискретизации интегральных уравнений 
Фредгольма первого рода \cite{ABT}. 

Когда  размерность $p$  велика, а собственные числа $\lambda_k$ достаточно быстро убывают очевидно, что наблюдения $\bar{Y}_k,\, k=1,\ldots,p$, (см. \eqref{eq.2}) могут содержать  существенную статистическую информацию о $\sigma^2$. Ясно, что эту  информацию можно извлечь, но при этом понятно также, что сделать это каким-то простым способом невозможно, 
поскольку  эти наблюдения содержат  неизвестные ``мешающие'' величины
$\sqrt{\lambda_k}\bar{\beta}_k$.

Понятно также, что для того, чтобы построить достаточно хорошие оценки $\sigma^2$,  надо вычесть из наблюдений $Y$  некоторую оценку вектора $X\beta$, а для этого нужно, естественно, оценивать $\beta$. 
Как мы видели  раньше, стандартную оценку максимального правдоподобия  применять не имеет смысла. Поэтому  далее в статье мы будем пользоваться  спектральными регуляризациями этой оценки, которые  вычисляются следующим образом: 
\begin{equation}\label{eq.4}
\widehat{\beta}_\alpha(Y)=H_\alpha(X^\top X)\bigl[(X^\top X)^{-1}X^\top Y\bigr];
\end{equation} 
здесь
\begin{itemize} 
\item
$  (X^\top X)^{-1}X^\top Y$ --  классическая оценка максимального правдоподобия неизвестного вектора $\beta$; 
\item $H_\alpha(\lambda)$ -- заданная функция $\mathbb{R}^+\rightarrow [0,1]$, индексированная параметром регуляризации $\alpha  \in \mathbb{R}^+$.
 \end{itemize} 

Матрица $H_\alpha(X^\top X)$, которая  сглаживает в спектральной области оценку максимального правдоподобия, формально определяется через спектральное разложение $X^\top X$ следующим образом:
\begin{equation*}
H_\alpha(X^\top X)=\sum_{k=1}^p H_\alpha(\lambda_k)\e_k \e_k^\top. 
\end{equation*}

Как правило, функции $H_\alpha(\cdot)$, используемые для регуляризации оценки максимального правдоподобия, обладают следующими естественными свойствами:
\begin{equation*}
\begin{split}
&0\le H_\alpha(\lambda)\le 1, \quad \text{при всех}\quad \lambda\ge 0,\ \alpha\ge 0;\\
&\lim_{\alpha\rightarrow 0}H_\alpha(\lambda)=1, \quad  \text{при любом}\quad\lambda>0;\\
&\lim_{\lambda\rightarrow 0}H_\alpha(\lambda)=0, \quad\text{при любом}\quad \alpha>0.
\end{split}
\end{equation*} 

Кроме того, семейство  функций $\{H_\alpha(\cdot),\, \alpha \in \mathbb{R}^+\}$ часто является упорядоченным, т.е. таким, что
для любых $\alpha,\alpha'\in \mathbb{R}^+$:\\ либо
\[
H_{\alpha}(\lambda)\le H_{\alpha'}(\lambda)\quad \text{при всех}\quad \lambda\in \mathbb{R}^+,
\]
либо 
\[
H_{\alpha'}(\lambda)\le H_{\alpha}(\lambda)\quad \text{при всех}\quad \lambda  \in \mathbb{R}^+.
\]
На исключительную полезность этого свойства при адаптивном (основанном на наблюдениях) выборе параметра регуляризации $\alpha$, по-видимому, впервые обратил внимание А. Кнайп \cite{K}.

На первый взгляд может показаться, что для вычисления   $\widehat{\beta}_\alpha(Y)$  нужно вычислять спектральное разложение $X^\top X$. Однако это не всегда так и на самом деле все зависит от того, каковы
функции $H_\alpha(\lambda)$.  Рассмотрим, например, метод регуляризации Тихонова 
\[
\widehat{\beta}^{\rm T}_\alpha(Y)=\arg\min_{\beta}\Bigl\{\|Y-X\beta\|^2+\alpha\|\beta\|^2\Bigr\}=(\alpha I+X^\top X)^{-1}X^\top Y.
\]  
Заметим, что 
\begin{itemize}
\item $\widehat{\beta}^{\rm T}_\alpha(Y)$ является решением линейного уравнения
\[
(\alpha I+X^\top X)\widehat{\beta}_{\rm T}(Y,\alpha)=X^\top Y;
\]
\item $\widehat{\beta}^{\rm T}_\alpha(Y)$ допускает следующее представление: 
\begin{equation*}
\begin{split}
\widehat{\beta}^{\rm T}_\alpha(Y)=&(\alpha I+X^\top X)^{-1}(X^\top X)(X^\top X)^{-1}X^\top Y\\ =& H_\alpha^{\rm T}(X^\top X)\bigl[(X^\top X)^{-1}X^\top Y\bigr],
\end{split}
\end{equation*}
где 
\[
H_\alpha^{\rm T}(\lambda)=\frac{\lambda}{\alpha +\lambda}.
\]
\end{itemize}
Таким образом, регуляризация Тихонова, с одной стороны, вычисляется как решение системы линейных уравнений, а с другой стороны,  представляет собой семейство упорядоченных спектральных регуляризаций. Примеры других   спектральных регуляризаций, которые вычисляются рекуррентно, можно найти, например, в \cite{EHN}. 

Отметим однако, что если  в качестве $H_\alpha(\lambda)$ используется, например, 
$H_\alpha(\lambda)=\mathbf{1}\{\lambda \ge \alpha \}$, то для того, чтобы построить соответствующие оценки из \eqref{eq.4}, необходимо уже вычислять спектральное разложение $X^\top X$. 

Далее, не оговаривая этого особо, будем предполагать, что для оценивания  $\beta$ используется семейство оценок
$\bigl\{\widehat{\beta}_\alpha(Y),\, \alpha\in\mathbb{R}^+\bigr\}$ из \eqref{eq.4} с упорядоченным семейством функций регуляризации $\{H_\alpha(\cdot),\, \alpha \in \mathbb{R}^+\}$.  Соответствующая 
$\widehat{\beta}_\alpha(Y)$   оценка $\sigma^2$
вычисляется как 
\begin{equation}\label{eq.5}
\begin{split}
\widehat{\sigma}^2_\alpha(Y)=&\frac{\bigl\|Y-X\widehat{\beta}_\alpha(Y)\bigr\|^2}{n}
\biggl[1-\frac{2{\rm tr}\bigl\{H_\alpha(X^\top X)\bigr\}-{\rm tr}\bigl\{[H_\alpha(X^\top X)]^2\bigr\}}{n}\biggr]^{-1}
.
\end{split}
\end{equation}

Таким образом, мы имеем семейство оценок дисперсии $\bigl\{\widehat{\sigma}^2_\alpha(Y),\, \alpha\in \mathbb{R}^+\bigr\}$ и, в принципе, задача состоит в том, чтобы выбрать в этом семействе наилучшую оценку на основе имеющихся наблюдений $Y$. 

Для этого  прежде всего нужно определить, что понимается под наилучшей оценкой. 
Представим себе 
идеальную ситуацию, когда $\beta=0$.  В этом случае мы наблюдаем чистый белый шум 
\begin{equation*}
Y_i=\sigma\xi_i, \quad i=1,\ldots,n
\end{equation*}  
и поэтому оценить  $\sigma^2$   в этом случае можно легко  с помощью  
\[
\widehat{\sigma}^2_\circ(\xi)\stackrel{\rm def}{=}\frac{\|\sigma\xi\|^2}{n}
=\sigma^2+\frac{\sigma^2}{n}\sum_{i=1}^n(\xi_i^2-1).\]
Понятно, что ничего принципиально лучшего, чем эта псевдо-оценка не существует. Поэтому мы будем измерять качество оценки $\widetilde{\sigma}^2(Y)$ величиной ее отклонения от 
$ \widehat{\sigma}^2_\circ(\xi)$, а именно, величиной
\[
\Delta\bigl(\widetilde{\sigma}^2\bigr)  ={n\bigl |\widetilde{\sigma}^2(Y)-\widehat{\sigma}^2_\circ(\xi)\bigr |}.
\] 
Очевидно, что чем меньше эта величина, тем лучше мы оцениваем уровень шума с помощью   $\widetilde{\sigma}^2(Y)$. Поэтому основная задача в этой работе -- минимизация $
\Delta\bigl(\widehat{\sigma}^2_\alpha\bigr)$ по $\alpha$ на основе имеющихся наблюдений. При этом будем предполагать, что параметр регуляризации $\alpha$ принадлежит отрезку $ \mathcal{A}=[\alpha_{ \min},\alpha_{\max}] $, а величина $\alpha_{\min}=\alpha_{\min}(n)>0$ зависит от 
$n$ и стремится к $0$ при $n\rightarrow\infty$.

Обозначим для краткости
\[
G_\alpha(\lambda)=2H_\alpha(\lambda)-[H_\alpha(\lambda)]^2
\quad 
\text{и} 
\quad
q_\alpha(n)=\biggl[1-\frac{1}{n}\sum_{i=1}^pG_\alpha(\lambda_i) \biggr]^{-1}-1.
\]

Для простоты  предположим, что
\begin{equation}\label{eq.6}
\frac{1}{n}\sum_{i=1}^p G_\alpha(\lambda_i)\ll 1.
\end{equation}

 Ясно, что без ограничения общности можно считать (см. \eqref{eq.2}), что $\lambda_k =0$ и $\bar{\beta}_k=0$ при $k=p+1,\ldots,n$. Тогда из \eqref{eq.2} и \eqref{eq.5} находим
\begin{equation}\label{eq.7}
\begin{split}
\widehat{\sigma}^2_\alpha(Y)=&
 \frac{1+q_\alpha(n)}{n}\sum_{i=1}^n\bigl[1-H_\alpha(\lambda_i)\bigr]^2\bar{Y}_i^2
\\ =& \frac{1+q_\alpha(n)}{n}\sum_{i=1}^n\bigl[1-H_\alpha(\lambda_i)\bigr]^2 \lambda_i\bar{\beta}_i^2
+ \frac{\sigma^2[1+q_\alpha(n)]}{n}\sum_{i=1}^n\bigl[1-G_\alpha(\lambda_i)\bigr] \xi_i'^2
\\ +&
 \frac{2\sigma[1+q_\alpha(n)]}{n}\sum_{i=1}^n\bigl[1-H_\alpha(\lambda_i)\bigr]^2 \xi_i' \sqrt{\lambda_i}\bar{\beta}_i
 .
\end{split}
\end{equation}
Воспользуемся также следующим простым тождеством:  
\begin{equation}\label{eq.8}
\begin{split}
 &\sigma^2[1+q_\alpha(n)]\sum_{i=1}^n\bigl[1-G_\alpha(\lambda_i)\bigr]\xi_i'^2 = \sigma^2\sum_{k=1}^n\xi'^2_k\\
&\qquad - \sigma^2[1+q_\alpha(n)]\sum_{k=1}^nG_\alpha(\lambda_k)(\xi_k'^2-1)
 +\sigma^2 q_\alpha(n)\sum_{k=1}^n (\xi_k'^2-1).
\end{split}
\end{equation}
Тогда с помощью \eqref{eq.7} и \eqref{eq.8}
 приходим к 
\begin{equation}\label{eq.9}
\begin{split}
\Delta\bigl(\widehat{\sigma}^2_\alpha\bigr) = &\biggl| [1+q_\alpha(n)]\sum_{i=1}^n\bigl[1-H_\alpha(\lambda_i)\bigr]^2 \lambda_i\bar{\beta}_i^2+ {\sigma^2}[1+q_\alpha(n)] \sum_{k=1}^nG_\alpha(\lambda_k)(1-\xi_k'^2)\\
+&\sigma^2q_\alpha(n) \sum_{i=1}^n (\xi_i'^2-1)+ 2\sigma[1+q_\alpha(n)]\sum_{i=1}^n\bigl[1-H_\alpha(\lambda_i)\bigr]^2 \xi_i' \sqrt{\lambda_i}\bar{\beta}_i\biggr|.
\end{split}
\end{equation}

Заметим, что второе слагаемое в правой части этого соотношения имеет 
порядок $\sqrt{\sum_{i=1}^nG_\alpha^2(\lambda_i)}$,  а третье  (см. \eqref{eq.6}) $n^{-1/2}\sum_{i=1}^nG_\alpha(\lambda_i)$. Поэтому   третьим слагаемым 
можно пренебречь если  в наряду с \eqref{eq.6}  выполняется также \\
\medskip
\noindent
\textbf{Условие A.} \textit {При всех} $\alpha \in \mathbb{R}^+$
\[
\sum_{i=1}^nG_\alpha(\lambda_i)\le K \sum_{i=1}^nG_\alpha^2(\lambda_i),
\]
\textit{где $K$ -- некоторая постоянная.}

Как мы увидим далее последнее слагаемое в правой части \eqref{eq.9} мало по сравнению с первым и поэтому из \eqref{eq.9}
получаем следующую аппроксимацию:
\begin{equation}\label{eq.10}
\begin{split}
\Delta\bigl(\widehat{\sigma}^2_\alpha\bigr) \approx &[1+q_\alpha(n)]\biggl| \sum_{i=1}^n\bigl[1-H_\alpha(\lambda_i)\bigr]^2 \lambda_i\bar{\beta}_i^2+ {\sigma^2} \sum_{k=1}^nG_\alpha(\lambda_k)(1-\xi_k'^2)\biggr|.
\end{split}
\end{equation}
 
Напомним, что наша цель -- выбрать параметр регуляризации $\alpha$  так, чтобы величина $\mathbf{E}\Delta\bigl(\widehat{\sigma}^2_\alpha\bigr)$ была бы  минимальной.  
Для этого нужно   на основе имеющихся наблюдений каким-то образом контролировать слагаемые в правой части \eqref{eq.10}. 

Рассмотрим сначала второе слагаемое, а именно, случайный процесс 
\[
\zeta(\alpha)=\sum_{k=1}^nG_\alpha(\lambda_k)(1-\xi_k'^2), \quad  \alpha\in \mathbb{R}^+.
\]
Нам этот процесс очевидно не известен поскольку у нас нет доступа к шумам $\xi_k', \, k=1,\ldots,n$. Единственное, что можно сделать  в этой ситуации -- это построить детерминированную верхнюю границу  для модуля этого процесса, т.е. найти в некотором смысле минимальную функцию 
$V(\alpha): (0,\alpha_{\max}]\rightarrow \mathbb{R}^+$ такую, что 
\begin{equation*}
\mathbf{E}\sup_{\alpha\le \alpha_{\max}}\bigl[|\zeta(\alpha)|-V(\alpha)\bigr]_+\le C \sqrt{\mathbf{E}\zeta^2(\alpha_{\max})},
\end{equation*} 
где здесь и далее $[x]_+=\max\{0,x\}$, а  $C$ обозначает постоянные  величины, значения которых могут меняться, но не зависят от параметров рассматриваемой задачи.  
Поиск такой функции -- нетривиальная задача, приближенное решение которой дает следующая теорема.
\begin{theorem}\label{th1} Пусть
\[
D(\alpha)= \sum_{k=1}^n G_\alpha^2(\lambda_k)
\]
и
\begin{equation} \label{eq.11}
{V}_\epsilon(\alpha)=(1+\epsilon) \sqrt{2D(\alpha)}\biggl\{
 \log \frac{D(\alpha)}{D(\alpha_{\max})}+2(1+\epsilon)\log\biggl[\frac{Q}{\epsilon^2}\log \frac{D(\alpha)}{D(\alpha_{\max})}\biggr]\biggr\}^{1/2},
\end{equation}
где
\begin{equation}\label{eq.12}
Q=\frac{4}{(\sqrt{2}-1)^2}.
\end{equation}
Тогда для любого  $\epsilon\in (0,1]$
\begin{equation}\label{eq.13}
\mathbf{E}\sup_{\alpha\le \alpha_{\max}}\bigl[|\zeta(\alpha)|-V_\epsilon(\alpha)\bigr]_+\le C\epsilon^{-1} \sqrt{D(\alpha_{\max})}.
\end{equation} 
\end{theorem}
Доказательство этого результата приведено в приложении. Оно по-сути основано на методе, который хорошо известен и используется при доказательстве закона повторного логарифма  \cite{Kol}.  

\bigskip
\noindent
\textbf{Замечание.} Функция $V_\epsilon(\alpha)$ из \eqref{eq.11} не является минимальной  детерминированной огибающей для $|\zeta(\alpha)|$. По-видимому, таковой является функция
\[
\widetilde{V}_\epsilon(\alpha)=
\sqrt{2D(\alpha)}\biggl\{\log \frac{D(\alpha)}{D(\alpha_{\max})}+ 2(1+\epsilon)\biggl[\log\log \frac{D(\alpha)}{D(\alpha_{\max})}+\log\frac{1}{\epsilon}\biggr]\biggr\}^{1/2},
\]
которая наряду с  $V_\epsilon(\alpha)$ обеспечивает выполнение
неравенства
\[
\mathbf{E}\sup_{\alpha\le \alpha_{\max}}\bigl[|\zeta(\alpha)|-\widetilde{V}_\epsilon(\alpha)\bigr]_+\le  C\epsilon^{-1} \sqrt{D(\alpha_{\max})}
\]
для любого $\epsilon\in (0,1]$.
К сожалению, строгого  доказательства этой гипотезы у нас нет.

\bigskip

Объединяя  \eqref{eq.13}  и \eqref{eq.10}, приходим к следующему  неравенству:
\begin{equation}\label{eq.14}
\begin{split}
\mathbf{E}\Delta\bigl(\widehat{\sigma}^2_{\widetilde\alpha}\bigr) \lesssim &\mathbf{E}[1+q_{\widetilde{\alpha}}(n)]\biggl\{\sum_{i=1}^n\bigl[1-H_{\widetilde{\alpha}}(\lambda_i)\bigr]^2 \lambda_i\bar{\beta}_i^2+\sigma^2V_\epsilon(\widetilde{\alpha})
\biggr\}+C\sigma^2\frac{\sqrt{ D(\alpha_{\max})}}{\epsilon},
\end{split}
\end{equation}
которое справедливо  для любого, зависящего от наблюдений, параметра регуляризации $\widetilde{\alpha}\le \alpha_{\max}$.

Очевидно, что мы хотели бы выбрать $\widetilde{\alpha}$ так, чтобы правая часть в \eqref{eq.14} была бы как можно меньше. Это идея приводит к  выбору
\begin{equation}\label{eq.15}
\widetilde\alpha(\beta)=\arg\min_{\alpha\in \mathcal{A}}[1+q_{\alpha}(n)] \biggl\{\sum_{i=1}^n\bigl[1-H_{{\alpha}}(\lambda_i)\bigr]^2 \lambda_i\bar{\beta}_i^2+\sigma^2{V}_\epsilon(\alpha)\biggr\}.
\end{equation}  

Очевидно, что этот параметр регуляризации нельзя использовать так как он зависит от неизвестного вектора $\beta$.   Поэтому наш следующий шаг состоит в том, чтобы оценить правую часть в \eqref{eq.15} по наблюдениям и тем самым  построить оценку для   $\alpha(\beta)$. Это можно сделать относительно просто,
заменив $ \lambda_i\bar{\beta}_i^2$ на $\bar{Y}_i^2-\sigma^2$ (см. \eqref{eq.2}). Таким образом,
мы получаем следующую оценку  для $\alpha(\beta)$: 
\begin{equation}\label{eq.16}
\begin{split}
\widetilde{\alpha}(Y)=\arg\min_{\alpha\in \mathcal{A}}\biggl\{[1+q_\alpha(n)] \biggl\{\sum_{i=1}^n\bigl[1-H_{\alpha}(\lambda_i)\bigr]^2 \bar{Y}_i^2+\sigma^2 {V}_\epsilon({\alpha})\biggl\}\\
+\sigma^2 [1+q_\alpha(n)] \sum_{i=1}^n G_\alpha(\lambda_i)-\sigma^2 n q_\alpha(n)\biggr\}.
\end{split}
\end{equation}

Эта оценка, как  видно, зависит от уровня шума, который мы и хотим оценить. Поэтому вместо $\sigma^2$ подставим в правую часть \eqref{eq.16}  ее простую оценку
\[
\frac{\bigl\|Y-X\widehat{\beta}_\alpha(Y)\bigr\|^2}{n}=\frac1n \sum_{i=1}^n\bigl[1-H_{\alpha}(\lambda_i)\bigr]^2 \bar{Y}_i^2=
\frac{\widehat{\sigma}^2_\alpha(Y)}{1+q_\alpha(n)}.
\]
Заметим также, что при выполнении \eqref{eq.6} 
\[
q_\alpha(n)\approx
\frac1n\sum_{i=1}^n G_\alpha(\lambda_i).
\]

  Таким образом,  \eqref{eq.16} приводит к следующему методу выбора параметра регуляризации:
\begin{equation}\label{eq.17}
\widehat{\alpha}(Y)=\arg\min_{\alpha\in \mathcal{A}} \biggl\{\widehat{\sigma}^2_\alpha(Y)
\biggl[1+\frac{{V}_\epsilon(\alpha)}{n}\biggr]\biggr\}.
\end{equation}  

Чтобы описать статистические свойства оценки $\widehat{\sigma}^2_{\widehat{\alpha}}(Y)$,  
нам потребуются следующие дополнительные обозначения: 
\begin{equation}\label{eq.18}
\begin{split}
&R_{\epsilon}(\alpha,\beta)\stackrel{\rm def}{=} \biggl[1+\frac{{V}_\epsilon(\alpha)}{n}\biggr]\biggl\{ [1+q_\alpha(n)] \sum_{i=1}^n\bigl[1-H_{{\alpha}}(\lambda_i)\bigr]^2 \lambda_i\bar{\beta}_i^2
+\sigma^2{V}_\epsilon(\alpha)\biggr\},
\\ &r_{\mathcal{A},\epsilon}(\beta)\stackrel{\rm def}{=}
\min_{\alpha\in\mathcal{A}}R_{\epsilon}(\alpha,\beta),\\ & 
\rho_{\mathcal{A},\epsilon}(\beta)\stackrel{\rm def}{=} \frac{\sigma^2\sqrt{D(\alpha_{\max})}}{r_{\mathcal{A},\epsilon}(\beta)}.
\end{split}
\end{equation}
\begin{theorem}\label{th-main}
Пусть  выполняется условие А и $\alpha_{\min},\, \alpha_{\max}$ таковы, что
\begin{equation}\label{eq.19}
\lim_{n\rightarrow\infty}\frac{D(\alpha_{\min})}{n}=0, \quad D(\alpha_{\max})\ge 5.
\end{equation} 
Тогда для любого $\gamma\in\bigl(0, \epsilon/ (1+\epsilon)\bigr)$, при всех $n\ge n_\gamma$ 
\begin{equation}\label{eq.20}
\begin{split}
\mathbf{E}\Delta(\widehat{\sigma}^2_{\widehat{\alpha}}) 
\le&  
\frac{ r_{\mathcal{A},\epsilon}(\beta)}{\gamma} \biggl\{1+  
\frac{1}{\sqrt[4]{D(\alpha_{\max})}} 
\log^{-1/8}\bigl[\rho_{\mathcal{A},\epsilon}(\beta)]
+\biggl[\frac{C}{(\epsilon-\gamma-\gamma\epsilon)\rho_{\mathcal{A},\epsilon}(\beta)}\biggr]^{1/2}\biggr\}^2.
\end{split}
\end{equation}
\end{theorem}

\bigskip
\noindent
\textbf{Замечания.} 
\begin{enumerate}
\item
Хорошо известно и легко проверить, что при $n\rightarrow\infty$
\[
\sqrt{n}\bigl[\sigma^2-\widehat{\sigma}^2_\circ(\xi)\bigr]\stackrel{\mathbf{D}}{\rightarrow}\sqrt{2} \sigma^2 \xi, 
\] 
где $\xi$ -- стандартная гауссовская случайная величина.
Поэтому неравенство \eqref{eq.20} фактически описывает члены второго порядка в разложении  $\widehat{\sigma}^2_{\widehat{\alpha}}(Y)$ по степеням $n^{-1/2}$. При этом выбор оптимального параметра регуляризации происходит с помощью минимизации этих членов.
 Отметим, что это типичная  ситуация  для задач семи-параметрического оценивания, см., например, \cite{GL,DGT}.  Подчеркнем также, что в отличии от этих статей, минимизация членов второго порядка в настоящей работе основана на имеющихся наблюдениях.  

Наиболее близкой по  математическим методам  к настоящей работе является статья \cite{LM}, хотя в ней рассматривается на первый совершенно другая статистическая задача оценивания квадратичного функционала. 

Отметим еще, что в современной статистике теория  оптимальности первого порядка семи-параметрических оценок хорошо разработана \cite{V}. В  ней  оптимальными оценками считаются все такие $\widetilde{\sigma}^{2}(Y)$, что 
\[
\sqrt{n}\big[\widetilde{\sigma}^{2}(Y)-{\sigma}^2]\stackrel{\mathbf{D}}{\rightarrow}\sqrt{2} \sigma^2 \xi, \quad n\rightarrow\infty.
\]   
 К сожалению, класс подобного рода оценок оказывается очень широким и в нем  невозможно определить наилучшую оценку или, что эквивалентно, сказать какой параметр регуляризации $\alpha$ будет наилучшим. 
\item
На первый взгляд кажется, что оптимальную оценку для $\sigma^2$ можно
получить если использовать наилучшую оценку для $X\beta$. На самом деле, как показывает теорема \ref{th-main}, это не так. Хорошо известно (см., например, \cite{K}), что наилучший параметр регуляризации $\alpha$ при оценивании $X\beta$ минимизирует 
 \[
\sum_{k=1}^n [1-H_\alpha(\lambda_k)]^2\lambda_k\bar{\beta}_k^2+\sigma^2\sum_{k=1}^n H_\alpha^2(\lambda_k),
\]
  в то время как при оценивании $\sigma^2$  оптимальное $\alpha$ минимизирует    
 \[
\sum_{k=1}^n [1-H_\alpha(\lambda_k)]^2 \lambda_k\bar{\beta}_k^2+\sigma^2 {V}_\epsilon(\alpha).
\]
Подчеркнем, что при  малых $\alpha$ (см. \eqref{eq.11})
$
{V}_\epsilon(\alpha)\ll \sum_{k=1}^n H_\alpha^2(\lambda_k).
$
\item Условие $\lim_{n\rightarrow\infty} D(\alpha_{\min})/n = 0$ является, по-видимому, техническим. На практике можно 
использовать $\alpha_{\min}=0$. К сожалению, математического доказательства этого предположения у нас нет.  
\item Величина $\rho_{\mathcal{A},\epsilon}(\beta)$, как правило, мала (см. параграф \ref{Re}). Поэтому, выражение в круглых скобках в правой части \eqref{eq.20} будет в этом случае близко к $1$.
\end{enumerate}

\section{Оценивание уровня шума в нелинейной регрессии \label{Re}}
В этом параграфе мы применим  теорему \ref{th-main}   к задаче минимаксного адаптивного оценивания $\sigma^2$
в нелинейной модели регрессии \eqref{eq.3}, предполагая, что число наблюдений $n$ велико. Термин ``минимаксный'' означает, что нас будет интересовать величина
\[
\bar{\Delta}( \widetilde{\sigma}^2,\mathcal{W}_2^m)=\max_{f\in \mathcal{W}_2^m}n\mathbf{E}\bigr|\widetilde\sigma^2(Y)-n^{-1}\sigma^2\|\xi\|^2\bigl|,
\] 
где $\widetilde\sigma^2(Y)$ -- некоторая оценка $\sigma^2$, построенная по наблюдениям \eqref{eq.3}.

Для того чтобы оценить неизвестную функцию $f(x), \, x\in [0,1]$, будем использовать семейство сглаживающих сплайнов 
\begin{equation}\label{eq.21}
\widehat{f}_\alpha(\cdot,Y)=\arg\min_{f}\biggl\{\frac{1}{n}\sum_{i=1}^n \bigl[Y_i-f(X_i)\bigr]^2+\alpha\int_0^1[f^{(m)}(x)]^2\, dx\biggr\}.
\end{equation}
 Поскольку этот метод является стандартным в непараметрической статистике, мы опустим его мотивацию и вычислительные аспекты,
которые  представлены, например, в \cite{GS}.

Хорошо известно, что получить простую  эквивалентную статистическую модель рассматриваемой
 задачи можно с помощью базиса Деммлера-Райнша \cite{DR}.  Базисные функции $\phi_{k}(x),\, 
k=1,\ldots,n,$ этого базиса обладают очень полезным свойством двойной ортогональности
\begin{equation*}
\begin{split}
&\frac{1}{n}\sum_{i=1}^n \phi_k(X_i)\phi_s(X_i)= \delta_{sk},\quad
\int_0^1 \phi_k^{(m)}(x)\phi_s^{(m)}(x)= \nu_k^n \delta_{ks},
\end{split}
\end{equation*}
где $\nu_1^n\le\nu_2^n\le \ldots\le \nu_n^n$ -- собственные числа базиса, которые зависят, естественно, от регрессоров $X_1,\ldots,X_n$. Чтобы упростить изложение, мы ограничимся асимптотическими равномерно распределенными регрессорами, т.е. такими, что
\[
\lim_{n\rightarrow\infty}\frac{1}{n}\sum_{i=1}^n \mathbf{1}\{X_i\le x\}=x,\quad x\in[0,1].
\]

Известно (см, например, \cite{S}),  что в этом случае и при  $n,k\rightarrow\infty$  
\begin{equation}\label{eq.22}
\nu_k^n=(1+o(1))(\pi k)^{2m}.
\end{equation}

  Легко проверить, что для эмпирических коэффициентов Фурье
\[
\bar{Y}_k=\frac{1}{n}\sum_{i=1}^n Y_i\phi_k(X_i)
\]  
справедливо представление
\begin{equation}\label{eq.23}
\bar{Y}_k= \bar{f}_k+\frac{\sigma}{\sqrt{n}}\xi_k', \quad k=1,\ldots,n,
\end{equation} 
где $\xi'$ -- стандартный белый гауссовский шум, а 
\[
\bar{f}_k=\frac{1}{n}\sum_{i=1}^n f(X_i)\phi_k(X_i).
\]

Заметим, что из свойства двойной ортогональности вытекает, что  сглаживающий сплайн из \eqref{eq.21} может быть представлен 
следующим образом:
\begin{equation*}
\widehat{f}_\alpha(x,Y)=\sum_{k=1}^n h_\alpha(\nu_k^n) \bar{Y}_k \phi_k(x),
\end{equation*}
где 
\[
h_\alpha(z)=\frac{1}{1+\alpha z}, \quad z>0,
\]
и что соболевский эллипсоид $\mathcal{W}_2^m$  в терминах  коэффициентов Фурье $\bar{f}_k,\, k=1,\ldots,n$, имеет вид
\begin{equation*}
\mathcal{W}_2^m=\biggl\{\bar{f}_k :\sum_{k=1}^n \nu_k^n \bar{f}_k^2\le L \biggr\}.
\end{equation*}

Поэтому статистическая модель наблюдений \eqref{eq.23} будет 
эквивалентна модели \eqref{eq.2} если  в последней положить
\begin{equation*}
\sigma=\frac{\sigma}{\sqrt{n}}, \quad \lambda_k =\frac{1}{\nu_k^n},\quad  \beta_k=\sqrt{\nu_k^n}\bar{f}_k.
\end{equation*}

Несложно проверить, воспользовавшись \eqref{eq.22}, что при $\alpha\rightarrow 0$
\begin{equation*}
\begin{split}
D(\alpha)=&(1+o(1))\sum_{k=1}^n \biggl[\frac{2}{1+\alpha(\pi k)^{2m}}-\frac{1}{[1+\alpha(\pi k)^{2m}]^2}\biggr]^2
=(1+o(1))\frac{K(m)\alpha^{-1/(2m)}}{\pi},
\end{split}
\end{equation*}
где 
\[
K(m)=\int_0^\infty\biggl[\frac{2}{1+x^{2m}}-\frac{1}{(1+x^{2m})^2}\biggr]^2\, dx
\]
и поэтому
\[
{V}_\epsilon(\alpha)=
(1+\epsilon+o(1))\alpha^{-1/(4m)}
\sqrt{\frac{K(m)}{\pi m}\log \frac{\alpha_{\max}}{\alpha}}.
\]

Оценка уровня шума, основанная на сглаживающих сплайнах, вычисляется как 
\[
\widehat{\sigma}^2_{\widehat{\alpha}}(Y)=\frac{1}{n}\sum_{i=1}^n [Y_i-\widehat{f}_{\widehat\alpha}(X_i,Y)]^2\biggl[1-\frac{W(\widehat\alpha)}{n}\biggr]^{-1},
\]
где
\begin{equation}\label{eq.24}
\widehat{\alpha}=\arg\min_{\alpha\in \mathcal{A}}\biggl\{
\sum_{i=1}^n [Y_i-\widehat{f}_\alpha(X_i,Y)]^2\biggl[1-\frac{W(\alpha)}{n}\biggr]^{-1}\biggl[1+
\frac{V_\epsilon(\alpha)}{n}\biggr]\biggr\},
\end{equation}
а
\[
W(\alpha)=\sum_{k=1}^n [2h_\alpha(\nu_k^2)-h_\alpha^2(\nu_k^2)].
\]

Для максимального  на классе $\mathcal{W}_2^m$ смещения этой оценки получаем следующую границу сверху:
\begin{equation*}
\sup_{f\in \mathcal{W}_2^m}\sum_{k=1}^n[1-h_\alpha(\nu_k^n)]^2\bar{f}_k^2=
L\max_{k}\frac{[1-h_\alpha(\nu_k^n)]^2}{\nu_k^n}=
L\alpha  \max_{k}\frac{\alpha\nu_k^n}{1+\alpha \nu_k^n} 
\le L\alpha.
\end{equation*}

Поэтому для любого $f\in \mathcal{W}_2^m$ при $n\rightarrow\infty$ (см. \eqref{eq.18})
\begin{equation*}
\begin{split}
r_{\mathcal{A},\epsilon}(f)\le&(1+\epsilon+o(1)) \min_{\alpha}\biggl\{L\alpha+\frac{\sigma^2}{n}\alpha^{-1/(4m)}\sqrt{\frac{K(m)}{\pi m}\log \frac{\alpha_{\max}}{\alpha}}\biggr\}\\
=&(1+\epsilon+o(1))\biggl(1+\frac{1}{4m}\biggr)\frac{\sigma^2}{n}\biggl(\frac{4mnL}{\sigma^2}\biggr)^{1/(4m+1)}\\
 \times&
\biggl[\frac{4K(m)}{\pi (4m+1)}\log \biggl(\frac{Ln}{\sigma^2}\biggr) \biggr]^{2m/(4m+1)}
\asymp \frac{\sigma^2}{n}\biggl(\frac{nL}{\sigma^2}\biggr)^{1/(4m+1)}
\biggl[\log \biggl(\frac{Ln}{\sigma^2}\biggr) \biggr]^{2m/(4m+1)}.
\end{split}
\end{equation*}
Отсюда сразу же вытекает (см. \eqref{eq.18}), что при $n\rightarrow\infty$
\[
\rho_{\mathcal{A},\epsilon}(f)\asymp\biggl(\frac{nL}{\sigma^2}\biggr)^{-1/(4m+1)}
\biggl[\log \biggl(\frac{Ln}{\sigma^2}\biggr) \biggr]^{-2m/(4m+1)}\rightarrow 0.
\]

Поэтому для оценки уровня шума из \eqref{eq.24} в силу теоремы \ref{th-main}  справедлива следующая асимптотическая (при $n\rightarrow\infty$)  верхняя граница:
\begin{equation}\label{eq.25}
\Delta(\widehat\sigma^2_{\widehat{\alpha}},\mathcal{W}_2^m)\stackrel{\rm asymp}{\le} \frac{C(m)}{\gamma} \frac{\sigma^2}{n}\biggl(\frac{nL}{\sigma^2}\biggr)^{1/(4m+1)}
\biggl[\log \biggl(\frac{Ln}{\sigma^2}\biggr) \biggr]^{2m/(4m+1)};
\end{equation}
здесь $\gamma<\epsilon/(1+\epsilon)$.

Несмотря на то, что рассматриваемая статистическая модель является стандартной в непараметричекой статистике и в ее рамках предложено много подходов к оцениванию $\sigma^2$, ничего не известно об оптимальности 
верхней границы \eqref{eq.25}. По-видимому, эту границу нельзя улучшить с точностью до множителя $(1+\epsilon)/\epsilon$, но доказательства этой гипотезы у нас, к сожалению, нет, хотя в качестве ее подтверждения  можно сослаться на \cite{EL}, где близкий факт доказан для задачи оценивания квадратичного функционала.

\section{Приложение}
\subsection{Вспомогательные результаты}

\begin{lemma}\label{lemma-0}
Пусть  $\xi_k, \, k=1,\ldots,n$, -- независимые $\mathcal{N}(0,1)$,
а $b_k,\, k=1,\ldots,n,$ -- детерминированная последовательность. Тогда для любого ${\widetilde\alpha}\in \mathbb{R}^+$, зависящего от
$\xi_k, \, k=1,\ldots,n,$ справедливо неравенство
\begin{equation*}
\mathbf{E}\biggl|\sum_{k=1}^n\bigl[1-H_{\widetilde\alpha}(\lambda_k)]^2 b_k \xi_s\bigg|
\le \biggl\{\mathbf{E}\sum_{k=1}^n\bigl[1-H_{\widetilde\alpha}(\lambda_k)]^4 b_s^2\biggr\}^{1/2}
.
\end{equation*}
\end{lemma}

\begin{lemma} \label{le1} 
Пусть 
\[
\zeta(\alpha_1,\alpha_2)=\sup_{\alpha\in [\alpha_1,\alpha_2]} \bigl[\zeta(\alpha)-\zeta(\alpha_2)\bigr].
\]
Тогда для любого $\lambda>0$  
\[
\mathbf{E}\exp\bigl[\lambda \zeta(\alpha_1,\alpha_2)\bigr]\le C\exp\bigl\{
Q\lambda^2[D(\alpha_1)-D(\alpha_2)]\bigr\},
\]
где постоянная $Q$ определена в \eqref{eq.12}.

\end{lemma}
Доказательство этих результатов вытекает из упорядоченности последовательностей $\bigl[1-H_\alpha(\lambda_k)]^2, \, k=1,\ldots,n $ и $G_\alpha(\lambda_k),
\, k=1,\ldots,n$,  (см., например, леммы 4 и 6 в  \cite{G16}). 

Нам также потребуется также еще один простой факт.
\begin{lemma}\label{l-entropy}
Пусть $\eta$ -- положительная случайная величина с ограниченным средним $m=\mathbf{E}\eta<\infty$, а $\theta$ -- случайная величина с ограниченным экспоненциальным моментом $\psi(\lambda)=\mathbf{E}\exp(\lambda\theta)<\infty $ для некоторого $\lambda>0$.
Тогда 
\begin{equation}\label{eq.26}
\mathbf{E}\eta\theta\le m \frac{H+\log[\psi(\lambda)]}{\lambda} ,
\end{equation}
где 
\[
H  = \mathbf{E}\frac{\eta}{m}\log\frac{\eta}{m}.
\]
\end{lemma}
\textsc{Доказательство.} Воспользуемся очевидным неравенством
\[
\theta\eta - z\exp(\lambda \theta)\le \max_x \bigl\{\eta x - z\exp(\lambda x)\bigr\}=
\frac{\eta}{\lambda}\log \frac{\eta}{z\lambda}-\frac{\eta}{\lambda},
\]
которое справедливо для любого $z>0$. Поэтому 
\begin{equation*}
\begin{split}
\mathbf{E} \theta\eta \le\min_{z>0}\biggl\{ z\psi(\lambda)-\frac{m}{\lambda}\log(z)+
\mathbf{E}\frac{\eta}{\lambda}\log \frac{\eta}{\lambda}-\frac{m}{\lambda}\biggr\}
= & \frac{mH+m\log[\psi(\lambda)]}{\lambda}.\quad \blacksquare
\end{split}
\end{equation*}

\begin{lemma}\label{lemma.4} Пусть выполнено условие А и \eqref{eq.19}. Тогда при $n\rightarrow\infty$
\begin{equation}\label{eq.27}
\begin{split}
\mathbf{E} \bigl\{q_{\alpha}(n)+ V_\epsilon(\alpha)\bigl[1+q_{\alpha}(n)\bigr]\bigr\}\biggl|\sum_{s=1}^n(\xi_s'^2-1)\biggr| \le& 
o(1)  \mathbf{E}V_\epsilon(\widehat\alpha)
\log^{-1/4}\frac{\mathbf{E}V_\epsilon(\widehat\alpha)}{\sqrt{D(\alpha_{\max})}}
.
\end{split}
\end{equation}
\end{lemma}
\textsc{Доказательство.}
Заметим, что если выполнено условие А и \eqref{eq.19}, то
\begin{equation}\label{eq.28}
\bigl\{q_{\alpha}(n)+ V_\epsilon(\alpha)\bigl[1+q_{\alpha}(n)\bigr]\bigr\} \le o(1) \frac{\sqrt{D(\widehat\alpha)}}{\sqrt{n}}.
\end{equation}
Далее воспользуемся леммой \ref{l-entropy}, положив 
\[
\theta=\sqrt{D(\widehat\alpha)}\quad \text{и} \quad \eta=\frac{1}{\sqrt{n}}
\sum_{i=1}^n (\xi_i'^2-1).
\] 
Поскольку в этом случае 
\[
\mathbf{E}\exp(\lambda\eta)=\exp\biggl[-\lambda\sqrt{n}-\frac{n}{2}\log\biggl(1-\frac{2\lambda}{\sqrt{n}}\biggr)\biggr],
\]
то при всех $\lambda \le \sqrt{n}/4$
\[
\mathbf{E}\exp(\lambda\eta)\le \exp(8\lambda^2).
\]
Поэтому из \eqref{eq.26} заключаем, что при  любом $\lambda \le \sqrt{n}/4$ выполнено неравенство
\begin{equation*}
\mathbf{E}\sqrt{D(\widehat\alpha)}\frac{1}{\sqrt{n}}
\sum_{i=1}^n (\xi_i'^2-1) \le \mathbf{E}\sqrt{D(\widehat\alpha)} \times
\biggl[\frac{1}{\lambda}\mathbf{E}\frac{\sqrt{D(\widehat\alpha)}}{\mathbf{E}\sqrt{D(\widehat\alpha)}}\log\frac{\sqrt{D(\widehat\alpha)}}{\mathbf{E}\sqrt{D(\widehat\alpha)}} +8\lambda\biggr].
\end{equation*}
Подставив в правую часть этого неравенства 
\[
\lambda =\frac{1}{2\sqrt{2}}\biggl[\mathbf{E}\frac{\sqrt{D(\widehat\alpha)}}{\mathbf{E}\sqrt{D(\widehat\alpha)}}\log\frac{\sqrt{D(\widehat\alpha)}}{\mathbf{E}\sqrt{D(\widehat\alpha)}}\biggr]^{1/2},
\]
приходим к 
\begin{equation*}\label{i1}
\mathbf{E}\sqrt{D(\widehat\alpha)}\frac{1}{\sqrt{n}}
\sum_{i=1}^n (\xi_i^2-1) \le 
C \sqrt{\mathbf{E}\sqrt{D(\widehat\alpha)}}
\biggl[\mathbf{E}{\sqrt{D(\widehat\alpha)}}\log\frac{\sqrt{D(\widehat\alpha)}}{\mathbf{E}\sqrt{D(\widehat\alpha)}}\biggr]^{1/2}.
\end{equation*}
Столь же просто проверяется, что
\begin{equation*}\label{i2}
\mathbf{E}\sqrt{D(\widehat\alpha)}\frac{1}{\sqrt{n}}
\sum_{i=1}^n (\xi_i^2-1) \ge 
-C \sqrt{\mathbf{E}\sqrt{D(\widehat\alpha)}}
\biggl[\mathbf{E}{\sqrt{D(\widehat\alpha)}}\log\frac{\sqrt{D(\widehat\alpha)}}{\mathbf{E}\sqrt{D(\widehat\alpha)}}\biggr]^{1/2}.
\end{equation*}
и, следовательно (см. \eqref{eq.28}),
\begin{equation}\label{eq.29}
\begin{split}
\mathbf{E} q_{\widehat\alpha}(n)\biggl|\sum_{s=1}^n(\xi_s'^2-1)\biggr| \le& 
o(1) \sqrt{\mathbf{E}D^{1/2}(\widehat\alpha)
\mathbf{E}V_\epsilon(\widehat\alpha)}
.
\end{split}
\end{equation}

 Нетрудно проверить также, что 
\begin{equation*}
D^{1/2}(\widehat\alpha)\le C V_\epsilon(\widehat\alpha)
\log^{-1/2}\frac{V_\epsilon(\widehat\alpha)}{\sqrt{D(\alpha_{\max})}} 
\end{equation*}
и поэтому  в силу неравенства Йенсена
\begin{equation*}
\mathbf{E}D^{1/2}(\widehat\alpha)\le C \mathbf{E}V_\epsilon(\widehat\alpha)
\log^{-1/2}\frac{\mathbf{E}V_\epsilon(\widehat\alpha)}{\sqrt{D(\alpha_{\max})}}. 
\end{equation*}
Таким образом
\begin{equation*}\label{equ.45}
\sqrt{\mathbf{E}D^{1/2}(\widehat\alpha)\mathbf{E} V_\epsilon(\widehat\alpha)}\le C \mathbf{E}V_\epsilon(\widehat\alpha)
\log^{-1/4}\frac{\mathbf{E}V_\epsilon(\widehat\alpha)}{\sqrt{D(\alpha_{\max})}}. 
\end{equation*}
Это неравенство и \eqref{eq.29} завершают доказательство леммы. $\quad\blacksquare$
\subsection{Доказательство теоремы \ref{th1}}

Заметим, что \eqref{eq.13} эквивалентно следующему неравенству:
\begin{equation*}
\mathbf{E}\sup_{\alpha\in \mathcal{A}}\bigl[\pm \zeta(\alpha)-V_\epsilon(\alpha)\bigr]_+\le C \frac{\sqrt{D(\alpha_{\max})}}{\epsilon}.
\end{equation*}

Выберем последовательность  $\alpha_k,\, k=0,1\ldots$, так, чтобы
\begin{equation}\label{eq.30}
D(\alpha_k)=(1+r)^k D(\alpha_{\max}),
\end{equation}
где величина $r>0$ будет определена позднее. 

Далее  воспользуемся следующими очевидными неравенствами:
\begin{equation}\label{eq.31}
\begin{split}
&\mathbf{E}\sup_{\alpha\in \mathcal{A}}\bigl[ \zeta(\alpha)-{V}_\epsilon(\alpha)\bigr]_+\le
\sum_{k=0}^\infty \mathbf{E} \sup_{\alpha\in [\alpha_{k+1},\alpha_{k}]}\bigl[ \zeta(\alpha)-{V}_\epsilon(\alpha)\bigr]_+\\ &\qquad \le
\sum_{k=0}^\infty \mathbf{E} \bigl[ \zeta(\alpha_k) 
-{V}_\epsilon(\alpha_k)+\sup_{\alpha\in [\alpha_{k+1},\alpha_{k}]}[\zeta(\alpha)-\zeta(\alpha_k)]\bigr]_+.
\end{split}
\end{equation}
и 
\begin{equation}\label{eq.32}
\begin{split}
&\mathbf{E}\sup_{\alpha\in \mathcal{A}}\bigl[- \zeta(\alpha)-{V}_\epsilon(\alpha)\bigr]_+\le
\sum_{k=0}^\infty \mathbf{E} \sup_{\alpha\in [\alpha_{k+1},\alpha_{k}]}\bigl[- \zeta(\alpha)-{V}_\epsilon(\alpha)\bigr]_+\\ &\qquad \le
\sum_{k=0}^\infty \mathbf{E} \bigl[- \zeta(\alpha_{k+1}) 
-{V}_\epsilon(\alpha_k)+\sup_{\alpha\in [\alpha_{k+1},\alpha_{k}]}[\zeta(\alpha_{k+1})-\zeta(\alpha)]\bigr]_+.
\end{split}
\end{equation}

Для того чтобы продолжить эти неравенства, заметим, что в силу неравенства Гельдера, для любых случайных величин
$\zeta_1, \zeta_2$ и любого $\lambda>0$ справедливы соотношения
\begin{equation}\label{eq.33}
\begin{split}
\mathbf{E}[\zeta_1-\zeta_2]_+\le &\lambda^{-1}\mathbf{E}\exp(-\lambda \zeta_2)\exp(\lambda \zeta_1)
\le \lambda^{-1}\bigl[\mathbf{E}\exp(-\lambda p\zeta_2)\bigr]^{1/p}\bigl[\mathbf{E}\exp(\lambda q\zeta_1)\bigr]^{1/q},
\end{split}
\end{equation}
где $p,q>1$ таковы,
что
\begin{equation}\label{eq.34}
\frac1p+\frac1q=1.
\end{equation}

В случае неравенства \eqref{eq.31} положим
\begin{equation*}
\begin{split}
\zeta_1=  \zeta(\alpha_k), \quad
\zeta_2=-{V}_\epsilon(\alpha_k)+\sup_{\alpha\in [\alpha_{k+1},\alpha_{k}]}[\zeta(\alpha)-\zeta(\alpha_k)]
.
\end{split}
\end{equation*}

Из леммы \ref{le1} (см. также \eqref{eq.30}) вытекает, что $\zeta_2$  суб-гауссовская случайная величина, точнее такая, для которой  неравенство 
\begin{equation}\label{eq.35}
\begin{split}
\mathbf{E}\exp(\lambda \zeta_2)\le& 
 \exp\bigl[-\lambda {V}_\epsilon(\alpha_k)
+Q r \lambda^2 D(\alpha_k) \bigr]
\end{split}
\end{equation}
выполняется для любого $\lambda>0$.

Далее легко проверить, что  $\zeta_1$ также суб-гауссовская случайная величина, т.е. для любого $\lambda>0$
\begin{equation}\label{eq.36}
\begin{split}
\mathbf{E}\exp(\lambda\zeta_1)\le& \exp\biggl\{\lambda
\sum_{s=1}^n G_{\alpha_k}(\lambda_s)-\frac12 \sum_{s=1}^n
\log\bigl[1+2\lambda G_{\alpha_k}(\lambda_s)\bigr]\biggr\}\\
\le& \exp\biggl\{\lambda^2
 \sum_{s=1}^n
\bigl[ G_{\alpha_k}(\lambda_s)\bigr]^2\biggr\}=\exp\bigl[\lambda^2 D(\alpha_k)\bigr].
\end{split}
\end{equation}

Следовательно, объединяя неравенства \eqref{eq.33}, \eqref{eq.35} и \eqref{eq.36},
находим
\begin{equation*}
\begin{split}
\mathbf{E} \bigl[ \zeta_1-\zeta_2\bigr]_+\le
C\lambda^{-1}\exp\bigl[-\lambda V_\epsilon(\alpha_k)
+(Q r  p+ q)\lambda^2 D(\alpha_k)\bigr]
\end{split} 
\end{equation*}
и подставив в это неравенство  
\[
\lambda=\frac{V_\epsilon(\alpha_k)}{2(Q r  p+ q) D(\alpha_k)},
\]
приходим к 
\begin{equation}\label{eq.37}
\begin{split}
\mathbf{E} \bigl[\zeta_1-\zeta_2\bigr]_+\le
\frac{2(Q r  p+ q)D(\alpha_k)}{V_\epsilon(\alpha_k)} \exp\biggl\{-\frac{V_\epsilon^2(\alpha_k)}{4(Q r  p+ q) D(\alpha_k)}\biggr\}.
\end{split} 
\end{equation}
Поскольку величина $q>1$ произвольна, выберем   
$
q=1+\sqrt{Qr }
$
и тогда из \eqref{eq.37} и \eqref{eq.34}
получим
\begin{equation}\label{eq.38}
\begin{split}
\mathbf{E} \bigl[\zeta_1-\zeta_2\bigr]_+\le C
\frac{C(1+\sqrt{Q r})^2D(\alpha_k)}{V_\epsilon(\alpha_k)} \exp\biggl\{-\frac{V_\epsilon^2(\alpha_k)}{4\bigr(1+\sqrt{Q r}\bigl)^2 D(\alpha_k)}\biggr\}.
\end{split} 
\end{equation}
Поэтому, если выбрать 
$r=\epsilon^2/Q$,
то в силу определения $V_\epsilon(\alpha)$ в \eqref{eq.11}, из \eqref{eq.38}
находим
\begin{equation*}
\begin{split}
\mathbf{E} \bigl[\zeta_1-\zeta_2\bigr]_+\le\frac{C \sqrt{D(\alpha_{\max})}}{k^{1+\epsilon}}
\end{split} 
\end{equation*}
и, таким образом, (см. \eqref{eq.31})
\[
\mathbf{E}\sup_{\alpha\in \mathcal{A}}\bigl[ \zeta(\alpha)-V_\epsilon(\alpha)\bigr]_+\le  C\epsilon^{-1}\sqrt{D(\alpha_{\max})}.
\]

Неравенство \eqref{eq.32}  продолжается практически аналогично. 
При этом  переопределим величины $\zeta_1, \, \zeta_2$ следующим образом:
\begin{equation*}
\begin{split}
\zeta_1= - \zeta(\alpha_{k+1}), \quad
\zeta_2=-V_\epsilon(\alpha_k)+\sup_{\alpha\in [\alpha_{k+1},\alpha_{k}]}\bigl[\zeta(\alpha_{k+1})-\zeta(\alpha)\bigr]
.
\end{split}
\end{equation*}

Случайная величина $\zeta_2$ как и ранее будет суб-гауссовской, 
но это уже не будет справедливо для $\zeta_1$, т.к.  
\begin{equation}\label{eq.39}
\begin{split}
\mathbf{E}\exp(\lambda \zeta_1)=
&
\mathbf{E}\exp\biggl[\lambda \sum_{s=1}^n G_{\alpha_{k+1}}(\lambda_s) (\xi_s'^2-1)\biggr]\\
 =& \exp\biggl\{- \lambda \sum_{s=1}^n G_{\alpha_{k+1}}(\lambda_s) 
 -\frac{1}{2}\sum_{s=1}^n \log\bigl[1-2\lambda G_{\alpha_{k+1}}(\lambda_s)  \bigr] \biggr\}.
\end{split}
\end{equation}
Однако,  это тождество  будет использоваться при 
 \begin{equation} \label{eq.40}
\begin{split}
\lambda =&\frac{V_\epsilon(\alpha_k)}{2(Q r  p+ q)D(\alpha_{k+1})}\le \frac{V_\epsilon(\alpha_k)}{2D(\alpha_{k+1})}\le \sqrt{\frac{\log[D(\alpha_k)/D(\alpha_{\max})]}{D(\alpha_k)}}\\=&\frac{1}{\sqrt{D(\alpha_{\max})}}\sqrt{\frac{\log[D(\alpha_k)/D(\alpha_{\max})]}{D(\alpha_k)/D(\alpha_{\max})}}\le \frac{1}{\sqrt{D(\alpha_{\max})}}.
\end{split}
\end{equation}
Поэтому, если $D(\alpha_{\max})>5$, то для оценки сверху правой части в \eqref{eq.39} можно использовать формулу Тейлора. Таким образом, для $\lambda$ из \eqref{eq.40}  получим
\begin{equation*}\label{t1.eqx}
\begin{split}
&\exp\biggl\{- \lambda \sum_{s=1}^n G_{\alpha_{k+1}}(\lambda_s) 
 -\frac{1}{2}\sum_{s=1}^n \log\bigl[1-2\lambda G_{\alpha_{k+1}}(\lambda_s)  \bigr] \biggr\}\\
&\quad \le \exp\biggl[\lambda^2 D(\alpha_{k+1})+C \frac{V_\epsilon^3(\alpha_k)}{D^2(\alpha_{k+1})}\biggr]
\le \exp\bigl[\lambda^2 D(\alpha_{k+1})+C\bigr].
\end{split}
\end{equation*}
Дальнейшее доказательство совершенно аналогично тому, что мы делали раньше и поэтому мы 
его опустим.

\subsection{Доказательство теоремы \ref{th-main}}
Обозначим для краткости
\begin{equation*}
\begin{split} 
&\widehat{\Delta}(\alpha,Y)=\biggl[1+\frac{V_\epsilon(\alpha)}{n}\biggr]\bigl[1+q_\alpha(n)\bigr]\bigl\|Y-X\widehat{\beta}_{\alpha}(Y)\bigr\|^2,\\ &
B(\alpha)=\sum_{s=1}^n\bigl[1-H_\alpha(\lambda_s)]^2\lambda_s\bar{\beta}_s^2 
\end{split}
\end{equation*} 
и (см. \eqref{eq.18})
\[
\alpha_\beta= \arg\min_{\alpha\in \mathcal{A}} R_{\epsilon}(\alpha,\beta) .
\]

Доказательство теоремы, по-существу, основано на  простом неравенстве
\begin{equation}\label{eq.41}
  \widehat{\Delta}(\widehat{\alpha},Y)\le \widehat{\Delta}(\alpha_\beta,Y)
\end{equation} 
и верхней границе, вытекающей из \eqref{eq.9} и теоремы 1 
\begin{equation*}\label{equ.38}
\begin{split}
\mathbf{E}\Delta\bigl(\widehat{\sigma}^2_{\widehat\alpha}\bigr) \le & R_\epsilon({\widehat\alpha},\beta)+  \frac{C\sigma^2\sqrt{D(\alpha_{\max})}}{\epsilon}\\
+&\sigma^2q_{\widehat\alpha}(n)\biggl| \sum_{i=1}^n (\xi_i'^2-1)\biggr|+ 2\sigma[1+q_{\widehat\alpha}(n)]\biggr|\sum_{i=1}^n\bigl[1-H_{\widehat\alpha}(\lambda_i)\bigr]^2 \xi_i' \sqrt{\lambda_i}\bar{\beta}_i\biggr|.
\end{split}
\end{equation*}

С помощью  лемм \ref{lemma-0} и \ref{lemma.4} это неравенство  можно продолжить следующим образом:
\begin{equation}\label{eq.42}
\begin{split}
\mathbf{E}\Delta\bigl(\widehat{\sigma}^2_{\widehat\alpha}\bigr) \le & R_\epsilon({\widehat\alpha},\beta)+  \frac{C\sigma^2\sqrt{D(\alpha_{\max})}}{\epsilon}\\
+&
C\sigma \sqrt{ \mathbf{E}B(\widehat{\alpha})} +
o(1)\sigma^2
 \mathbf{E}V_\epsilon(\widehat\alpha)
\log^{-1/4}\frac{\mathbf{E}V_\epsilon(\widehat\alpha)}{\sqrt{D(\alpha_{\max})}}
.
\end{split}
\end{equation}

Чтобы воспользоваться  этим неравенством, нужно выяснить как связаны между собой случайный процесс  $\widehat{\Delta}(\alpha,Y),\, \alpha\in\mathcal{A},$  и  функция $R_{\epsilon}(\alpha,\beta),\, \alpha\in\mathcal{A}$.
Эту связь определяет следующее тождество:
\begin{equation}\label{eq.43}
\begin{split}
\widehat{\Delta}(\alpha,Y)-\sigma^2\|\xi\|^2
=&\biggl[1+\frac{{V}_\epsilon(\alpha)}{n}\biggr] \bigl[1+q_\alpha(n)\bigr]B(\alpha)  \\ +&2\sigma \bigl[1+q_\alpha(n)\bigr]\biggl[1+\frac{{V}_\epsilon(\alpha)}{n}\biggr]\sum_{s=1}^n\bigl[1-H_\alpha(\lambda_s)]^2\sqrt{\lambda_s}\bar{\beta}_s \xi_s'\\
+&\sigma^2 \biggl[1+\frac{{V}_\epsilon(\alpha)}{n}\biggr]\biggl[ V_\epsilon(\alpha)-\sum_{s=1}^nG_\alpha(\lambda_s)(\xi_s'^2-1)\biggr]
\\ +&\sigma^2\bigl\{q_{\alpha}(n)+ V_\epsilon(\alpha)\bigl[1+q_{\alpha}(n)\bigr]\bigr\}\sum_{s=1}^n(\xi_s'^2-1).
\end{split}
\end{equation}

Отсюда  и из \eqref{eq.41} сразу же получаем
\begin{equation}\label{eq.44}
\begin{split}
\mathbf{E}\bigl[\widehat{\Delta}(\widehat\alpha,Y)-\sigma^2\|\xi\|^2\bigr]
\le \mathbf{E}\bigl[\widehat{\Delta}(\alpha_\beta,Y)-\sigma^2\|\xi\|^2\bigr]
= R_{\epsilon}(\alpha_\beta,\beta)=r_{\mathcal{A},\epsilon}(\beta).
\end{split}
\end{equation}

Нам потребуется   нижняя граница для левой части в этом неравенстве. 
Ее можно получить подставив в \eqref{eq.43} 
\[
V_\epsilon(\alpha)=\bigl(1-\gamma\bigr)
V_\epsilon(\alpha)+\gamma V_\epsilon(\alpha);
\]
здесь $\gamma\in (0,\epsilon/(1+\epsilon))$. 
Заметим  (см. \eqref{eq.11}), что 
\[
\bigl(1-\gamma \bigr)V_\epsilon(\alpha) \ge V_{\epsilon -\gamma-\gamma\epsilon}(\alpha).
\]
При выполнении  условия А и \eqref{eq.19} очевидно, что $q_\alpha(n)< 1$. Поэтому  из  \eqref{eq.43} с помощью  \eqref{eq.29}, лемм \ref{lemma-0}, \eqref{lemma.4} и теоремы \ref{th1} приходим к неравенству 
\begin{equation*}
\begin{split}
&\mathbf{E}\bigl[\widehat{\Delta}(\widehat\alpha,Y)-\sigma^2\|\xi\|^2\bigr]
\ge \mathbf{E} \biggl[1+\frac{{V}_\epsilon(\widehat\alpha)}{n}\biggr] \bigl\{
\bigl[1+q_{\widehat{\alpha}}(n)\bigr]B(\widehat{\alpha}) +\sigma^2 \gamma  V_\epsilon(\widehat\alpha) \bigr\}\\
&\qquad -C\sigma\sqrt{\mathbf{E}B(\widehat{\alpha})}
- C \sigma^2\frac{\sqrt{D(\alpha_{\max})}}{(\epsilon-\gamma-\gamma\epsilon)_+}
- 
o(1)\sigma^2\mathbf{E}V_\epsilon(\widehat\alpha)
\log^{-1/4}\frac{\mathbf{E}V_\epsilon(\widehat\alpha)}{\sqrt{D(\alpha_{\max})}}.
\end{split}
\end{equation*}

Воспользовавшись этим неравенством и \eqref{eq.44}, получим
\begin{equation}\label{eq.45}
\begin{split}
r_{\mathcal{A},\epsilon}(\beta)
\ge& \sigma^2\gamma  \mathbf{E}V_\epsilon(\widehat\alpha)   -  C \sigma^2\frac{\sqrt{D(\alpha_{\max})}}{(\epsilon-\gamma-\epsilon\gamma)_+} -o(1)\sigma^2\mathbf{E}V_\epsilon(\widehat\alpha)
\log^{-1/4}\frac{\mathbf{E}V_\epsilon(\widehat\alpha)}{\sqrt{D(\alpha_{\max})}}
\end{split}
\end{equation}
и, решая это неравенство относительно  $\mathbf{E}V_\epsilon(\widehat\alpha)$, приходим к
\begin{equation*}
\begin{split}
\mathbf{E}\frac{V_\epsilon(\widehat\alpha)}{\sqrt{D(\alpha_{\max})}} &\le (1+o(1))\biggl[\frac{ r_{\mathcal{A},\epsilon}(\beta)}{\sigma^2\gamma\sqrt{D(\alpha_{\max})}} +  \frac{C}{\gamma(\epsilon-\gamma-\epsilon\gamma)_+}\biggr].
\end{split}
\end{equation*}

С  помощью  этого неравенства получим
\begin{equation}\label{eq.46}
\begin{split}
\mathbf{E}V_\epsilon(\widehat\alpha)
\log^{-1/4}\frac{\mathbf{E}V_\epsilon(\widehat\alpha)}{\sqrt{D(\alpha_{\max})}}\le&
\frac{Cr_{\mathcal{A},\epsilon}(\beta)}{\gamma\sigma^2\sqrt{D(\alpha_{\max})}} \biggl[1
 + \frac{C\sigma^2\sqrt{D(\alpha_{\max})}}{(\epsilon-\gamma-\epsilon\gamma)_+ r_{\mathcal{A},\epsilon}(\beta) }\biggr]\\  \times& 
\log^{-1/4}\biggl(\frac{ r_{\mathcal{A},\epsilon}(\beta)}{\sigma^2\sqrt{D(\alpha_{\max})}}\biggr).
\end{split}
\end{equation}

Подставив в \eqref{eq.45}  неравество  \eqref{eq.46}  получаем, что при достаточно больших $n$ 
\begin{equation*}
\begin{split}
r_{\mathcal{A},\epsilon}(\beta)
\ge& \mathbf{E}
\bigl[1+q_{\widehat{\alpha}}(n)\bigr]\bigl[B(\widehat{\alpha})
 +\gamma \sigma^2 V_\epsilon(\widehat\alpha)\bigr] \\
-&C\sigma\Bigl\{\mathbf{E}\bigl[1+q_{\widehat{\alpha}}(n)\bigr]
\bigl[B(\widehat{\alpha}) +\gamma\epsilon\sigma^2 V_\epsilon(\widehat\alpha)\bigr]\Bigr\}^{1/2}
- \frac{C \sigma^2\sqrt{D(\alpha_{\max})}}{(\epsilon-\gamma-\gamma\epsilon)_+}\\ -&o(1) 
\frac{r_{\mathcal{A},\epsilon}(\beta)}{\gamma\sqrt{D(\alpha_{\max})}} 
\log^{-1/4}\biggl[\frac{r_{\mathcal{A},\epsilon}(\beta)}{\sigma^2\sqrt{D(\alpha_{\max})}}  
\biggr]
.
\end{split}
\end{equation*}
Отсюда с помощью простых преобразований приходим к неравенству
\begin{equation*}
\begin{split}
&\mathbf{E}\bigl[1+q_{\widehat{\alpha}}(n)\bigr]
\bigl[B(\widehat{\alpha}) +\gamma \sigma^2 V_\epsilon(\widehat\alpha)\bigr]
\\ &\qquad 
\le r_{\mathcal{A},\epsilon}(\beta)\biggl\{1+  
\frac{o(1)}{\sqrt[4]{D(\alpha_{\max})}\sqrt{\gamma}} 
\log^{-1/8}\bigl[\rho_{\mathcal{A},\epsilon}(\beta)]
+\biggl[\frac{C}{(\epsilon-\gamma-\gamma\epsilon)_+\rho_{\mathcal{A},\epsilon}(\beta)}\biggr]^{1/2}\biggr\}^2.
\end{split}
\end{equation*}
и, таким образом,  
\begin{equation*}
\begin{split}
\mathbf{E}R_\epsilon(\widehat{\alpha},\beta) 
\le& \frac{r_{\mathcal{A},\epsilon}(\beta)}{\gamma}\biggl\{1+  
\frac{o(1)}{\sqrt[4]{D(\alpha_{\max})}\sqrt{\gamma}} 
\log^{-1/8}\bigl[\rho_{\mathcal{A},\epsilon}(\beta)]
+\biggl[\frac{C}{(\epsilon-\gamma-\gamma\epsilon)_+\rho_{\mathcal{A},\epsilon}(\beta)}\biggr]^{1/2}\biggr\}^2.
\end{split}
\end{equation*}

Наконец, подставляя это неравенство и \eqref{eq.46} в \eqref{eq.42}, с помощью несложных вычислений приходим  
к \eqref{eq.20}.

\begin{center}{\large  СПИСОК ЛИТЕРАТУРЫ} \end{center}


\vspace{5mm}

\begin{flushleft}
{\small {\it Голубев Георгий Ксенофонтович} \\
CNRS, Aix-Marseille Universit\'e, I2M,\\
Институт проблем передачи информации 
им. А.А. Харкевича РАН \\
 {\tt golubev.yuri@gmail.com}}

\bigskip 

{\small {\it Крымова Екатерина Александровна} \\
Институт проблем передачи информации
 им. А.А. Харкевича РАН, \\
 Duisburg-Essen University,\\
 {\tt ekkrym@gmail.com}} 
\end{flushleft}%


\begin{enumerate}
\bibitem{ABT}
{\it Aster, R., Borcher, B., Thurber, C. } Parameter Estimation and Inverse Problems.  Elsevier 2013.
\bibitem{K}
\textit{Kneip,~A.} 
Ordered linear smoothers //
Annals of Statist.  1994. V.  22.  № 2.  P. 835--866.


\bibitem{EHN}  \textit{Engl, H.W.,  Hanke, M., and  Neubauer, A.} 
 Regularization of Inverse Problems // Mathematics and its Applications. V. 375.
Dordrecht: Kluwer Academic Publishers Group 1966.


\bibitem{Kol}
 {\it  Kolmogoroff A.} Uber das Gesetz des iterierten Logarithmus // Mathematische Annalen. 1929. V. 101. P. 126--135.

\bibitem{GL}  {\it Golubev G., Levit B.} On the second order minimax estimation of
distribution function // Math. Methods Statist. 1996. V. 5. P.
1--31.

\bibitem{DGT}{\it Dalalyan A, Golubev G. and Tsybakov A.} Bayesian maximum
likelihood and  semi\-parametric second order efficiency // {Ann.
Statist. } 2006. V. 34, № 1. P. 169--201.

\bibitem{LM} {\it Laurent B. and Massart P} Adaptive estimation of a quadratic functional by model selection // The Annals of Statist. 2000. V. 28. № 5. P. 1302--1338. 





 
\bibitem{V} {\it Van der Vaart A.} Simiparametric statistics // 
 Lectures on Probability Theory and Statistics. Ecole d'et\'e de probalit\'es de Saint-Flour XXIX-1999 P. 312--425.

\bibitem{GS}
\textit{Green, P.\,J. and  Silverman, B.\,W.}
 Nonparametric Regression and Generalized Linear Models.
A roughness penalty approach,  Chapman and Hall, 1994.


\bibitem{DR}
{\it Demmler A., Reinsch C.} Oscillation matrices with spline smoothing // Numerische Mathematik 1975.  V. 24. № 5. P. 375–-382.
\bibitem{S}
{\it Speckman P.}
Spline smoothing and optimal rates of convergence in nonparametric regression models. //
Annals of Statist. 1985. V. 13 № 3. P. 970--983.
\bibitem{G16}{\it Голубев Г. К.} Концентрации рисков выпуклых комбинаций линейных оценок. // Проблемы перед. информ. 2016. Т. 52. № 4. С. 31--48.

\bibitem{EL} {\it Efromovich S. and Low M.} On optimal adaptive estimation of a quadratic functional. 
// The Annals of Statist. 1996. V. 24. № 3. P. 1106--1125. 


\end{enumerate}
\end{document}